\documentclass[11pt]{article}  

\usepackage[dvips]{epsfig}
\usepackage{latexsym}
\usepackage{amsfonts,amsmath,amssymb}
\newtheorem{corollary}{Corollary}
\newtheorem{lemma}{Lemma}

\newtheorem{theorem}{Theorem}
\newtheorem{proposition}{Proposition}

\newcommand{\qed}{\mbox{$\Diamond$}\vspace{\baselineskip}}
\newenvironment{proof}{\noindent{\bf Proof:}}{\qed}
\newcommand{\Var}{\hbox {Var}}

\begin{document}
\title{Real Zeros and Normal Distribution for statistics
on Stirling permutations defined by Gessel and Stanley}

\author{Mikl\'os B\'ona \\
Department of Mathematics \\
University of Florida\\
 Gainesville FL 32611-8105\\
bona@math.ufl.edu \thanks{Partially supported
by an NSA Young Investigator Award.}}

\date{}

\maketitle

\begin{abstract} We study Stirling permutations defined by Gessel and
Stanley in \cite{stangess}. We prove that their generating function according
to the number of descents has real roots only. We use that fact to prove
that the distribution of these descents, and other, equidistributed 
statistics on these objects converge to a normal distribution.
\end{abstract}

\section{Introduction}
In \cite{stangess} Ira Gessel and Richard Stanley defined an interesting
class of multiset permutations called {\em Stirling Permutations}.
 Let $Q_k$ denote the set of all permutations
of the multiset $\{1,1,2,2,\cdots ,n,n\}$ in which for all $i$, 
 all entries between the two occurrences of $i$ are larger than $i$. 
For instance, $Q_2$ has three elements, namely 1122, 1221, and 2211. 
It is not difficult to see that $Q_n$ has $1\cdot 3\cdot \cdots \cdot
(2n-1)=(2n-1)!!$ elements. 
Gessel and Stanley then proved many enumerative results for these permutations
and showed several connections between these and
 other combinatorial objects, such as set partitions. 

Counting Stirling permutations by descents, the authors
of \cite{stangess} found a recurrence relation similar to the recurrence
relation  known for classic permutations. In this paper, we will continue in
that direction. First, we show the simple but interesting fact that 
on $Q_n$ the descent and the {\em plateau} statistics, to be defined in the
next section, are equidistributed.
Then we prove that for any fixed $n$, the generating polynomial of all
Stirling permutations in $Q_n$ with respect to the descent statistic has real
roots only. This is analogous to the well-known case
(see  Theorem 1.33 of \cite{bona})
 of classic permutations, namely the
result that all the roots of  Eulerian polynomials are real.  Finally, 
we apply a classic result of Bender to use this real roots property to prove
that the descents of Stirling permutations in $Q_n$ are normally distributed.

\section{Stirling Permutations and Real Zeros}
Let $q=a_1a_2\cdots a_{2n}\in Q_n$ be a Stirling permutation. 
Let the index $i$ be called an {\em ascent} of $q$ 
if $i=0$ or $a_i<a_{i+1}$, let $i$
be called a {\em descent} of $q$
 if $i=2n$ or $a_i>a_{i+1}$, and let $i$ be called a
{\em plateau} of $q$
 if $a_i=a_{i+1}$. It is obvious that the ascent and descent
statistics are equidistributed, since reversing an element of $Q_n$ turns
ascents into descents and vice versa. It is somewhat less obvious that the
plateau statistic is also equidistributed with the previous two. This fact, 
and a reason for it, are 
the content of the next proposition. Note that its first identity, 
(\ref{staneq}), was proved in \cite{stangess}. 

\begin{proposition} \label{stanp}
Let $C_{n,i}$ be the number of elements of $Q_n$ with $i$ descents.
Then for all positive integers $n,i\geq 2$, we have
\begin{equation} \label{staneq} C_{n,i}=iC_{n-1,i}+(2n-i)C_{n-1,i-1}.
\end{equation}
Similarly, let $c_{n,i}$ be number of elements of $Q_n$ with $i$
plateaux. Then  for all positive integers $n,i\geq 2$, we have
\begin{equation} \label{plateq} c_{n,i}=ic_{n-1,i}+(2n-i)c_{n-1,i-1}.
\end{equation}
In particular, since $C_{1,1}=c_{1,1}=1$ and $C_{1,0}=c_{1,0}=0$, the 
identity \begin{equation} \label{allthesame}
C_{n,i}=c_{n,i} \end{equation} holds. 
\end{proposition}

\begin{proof}
There are two ways to obtain an element of $Q_n$ from an element $p\in 
Q_{n-1}$ by inserting two copies of $n$ into consecutive positions. Either
$p$ must have $i$ descents, and then we insert the two copies of $n$ into a
descent, or $p$ has $i-1$ descents, and then we insert the two consecutive
copies of $n$ into one of the $(2n-1)-(i-1)=2n-i$ positions that are not
descents.  

The argument proving (\ref{plateq}) is analogous.
\end{proof}

\begin{corollary} \label{expect} 
On average, elements of $Q_n$ have
$(2n+1)/3$ ascents, $(2n+1)/3$ descents, and $(2n+1)/3$ plateaux. 
\end{corollary}

Proposition \ref{stanp} enables us to prove a strong result on
 the roots of the polynomials $\sum_{i=1}^n C_{n,i}x^i$. The method we use
follows an idea of  H. Wilf  (\cite{hwilf}, \cite{bona} Theorem 1.33) who
used it on classic permutations. 

\begin{theorem} \label{realzeros} 
Let $C_{n}(x)=\sum_{i=1}^n C_{n,i}x^i$. Then for all positive
integers $n$, the roots of the polynomial $C_{n}(x)$ are all real, distinct,
and non-positive. 
\end{theorem}

\begin{proof} For $n=1$, one sees that $C_1(x)=x$, and the statement holds. 
For $n=2$, one sees that $C_2(x)=2x^2+x=x(2x+1)$, and so the statement again
holds.

For $n\geq 3$, recurrence relation (\ref{staneq}) implies 
\begin{equation} \label{first} 
C_{n}(x)=(x-x^2)C_{n-1}'(x)+(2n-1)xC_{n-1}(x)\end{equation}
as can be seen by equating coefficients of $x^i$. The right-hand
side is similar to the derivative of a product, which suggests the following
rearrangement
\begin{equation} \label{wilf} 
C_{n}(x)=x(1-x)^{2n}\frac{d}{dx}\left ((1-x)^{1-2n}C_{n-1}(x) \right).
\end{equation} 
Let us now assume inductively that the roots of $C_{n-1}(x)$ are real, 
distinct and non-positive. Clearly, $C_n(x))$ vanishes at $x=0$. Furthermore,
by Rolle's theorem, (\ref{wilf}) shows that $C_n(x)$ has a root between 
any pair of consecutive roots of $C_{n-1}(x)$. This counts for $n-1$ roots
of $C_n(x)$. So the last root must also be real, since complex roots of 
polynomials with real coefficients must come in conjugate pairs.  

There remains to show  that the
last root of $C_{n}(x)$ must be on the right of the rightmost root of
$C_{n-1}$. Consider (\ref{first}) at the rightmost root $x_0$
 of $C_{n-1}$.
As $x_0$ is negative, we know that $x_0-x_0^2<0$, and so 
$C_n(x_0)$ and $C_{n-1}'(x_0)$ have opposite signs. The claim now follows, 
since in $-\infty$, the polynomials  $C_n(x)$ and $C_{n-1}'(x)$ must
converge to the same (infinite) limit as their degrees are of the same
parity. As $C_{n-1}'(x)$ has no more roots on the right of $x_0$, the 
polynomial $C_{n}(x)$ must have one.
\end{proof}

Note that we have in fact proved that the roots of $C_{n-1}(x)$ and
$C_n(x)$ are interlacing, so the sequence $C_1,C_2, \cdots $ is a {\em
Sturm sequence.} 

As an immediate application of the real zeros property, we can determine
where peak (or peaks) of the sequence $C_{n,1}, C_{n,2},\cdots, C_{n,n}$ is.
Our tool in doing so is the following theorem of Darroch.

\begin{theorem} \label{darroch} \cite{darroch}
Let $A(x)=\sum_{k=0}^n a_kx^k$ be a polynomial that has
real roots only that satisfies $A(1)>0$. Let $m$ be an index so that
$a_m=\max_{0\leq i\leq n} a_i$.
 Let $\mu=A'(1)/A(1)$.
Then we have
\[|\mu -m| <1 .\] 
\end{theorem}

In particular, a sequence with the real zeros property can have at most
two peaks. Note that $A'(1)=\sum_{i= 0}^n ia_i$ and $A(1)=\sum_{i=0}^n a_i$,
therefore $A'(1)/A(1)$ is nothing else but the weighted average of the
coefficients $a_i$, with $i$ being the weight of $a_i$. So in the particular
case when $A(x)=C_n(x)$, we have
\begin{eqnarray*} \frac{C_n'(1)}{C_n(1)} & = 
& \frac{\sum_i iC_{n,i}}{\sum_{i} C_{n,i}} \\
& = & \sum_i i \cdot \frac{C_{n,i}}{(2n-1)!!} \\
& = & \frac{2n+1}{3},
\end{eqnarray*} 
where the last step follows from Corollary \ref{expect}. Indeed,  
$\frac{C_{n,i}}{(2n-1)!!}$ is just the probability that a randomly selected
Stirling permutation of length $n$ has exactly $i$ descents, so 
$sum_i i \cdot \frac{C_{n,i}}{(2n-1)!!}$ is just the expected number of
descents in such permutations. 

Therefore, by Theorem \ref{darroch}, we obtain the following result.
\begin{theorem}
Let $i$ be an index so that $C_{n,i}=\max_k C_{n,k}$. Then
\begin{enumerate}
\item $i=(2n+1)/3$ if $(2n+1)/3$ is an integer, and
\item $i=\lfloor (2n+1)/3 \rfloor$ or  $i=\lceil (2n+1)/3 \rceil$ if
 $(2n+1)/3$ is not an integer.
\end{enumerate}
\end{theorem}

\section{Stirling Permutations and Normal Distribution}
In this section, we prove that the plateaux (equivalently ascents, 
equivalently, descents) of Stirling permutations are normally distributed.
Our main tool is the following result of Bender. 
Let $X_n$ be a random variable, and let $a_n(k)$ be
a triangular array of non-negative real numbers, $n=1, 2,\cdots $, and
$1\leq k\leq m(n)$ so that
\[P(X_n=k)=p_n(k)=\frac{a_n(k)}{\sum_{i=1}^{m(n)} a_n(i)}.\]
Set $g_n(x)=\sum_{k=1}^{m(n)} p_n(k)x^k$.

We need to introduce some notation for transforms of the random variable
$Z$. Let $\bar{Z}=Z-E(Z)$, let $\tilde{Z}=\bar{Z}/\sqrt{\Var( Z)}$, and let
$Z_n\rightarrow N(0,1)$ mean that $Z_n$ converges in distribution to the 
standard normal variable.

\begin{theorem} \label{bender} \cite{bender}
Let $X_n$ and $g_n(x)$ be as above. If $g_n(x)$ has real roots only, and
\[\sigma_n =\sqrt{\Var( X_n)}\rightarrow \infty,\] then
 $\tilde{X}_n\rightarrow N(0,1)$.
\end{theorem}

See \cite{cannorm} for related results.

 We want 
to use Theorem \ref{bender} to prove that the plateaux of permutations
in $Q_n$ are normally distributed. Because of Theorem \ref{realzeros}, all
we need for that is to prove that the variance of the number of these
plateaux converges to infinity as $n$ goes to infinity. We will accomplish
more by proving an explicit formula for this variance. 
In order to state that formula, 
let $Y_{n,i}$ be the indicator random variable of the event that in a
randomly selected element of $Q_n$, the two copies of $i$ are consecutive,
that is, they form a plateau. Note that $P(Y_{n,n}=1)=E(Y_{n,n})=1$. 
Set $Y_n=\sum_{i=1}^n Y_{n,i}$. 

\begin{theorem} \label{varform} For all positive integers $n$, the 
equality
\begin{equation} \label{explicit} \Var(Y_n)=\frac{2n^2-2}{18n-9}
\end{equation} holds.
\end{theorem}

\begin{proof} We are going to use the identity $\Var(Y_n)=E(Y_n^2)-E(Y_n)^2$.
We have seen in Corollary \ref{expect} that $E(Y_n)=\frac{2n+1}{3}$. 
Let $s_n=E(Y_n^2)$. The key element of our computations is the following
lemma.

\begin{lemma} \label{recu} For all positive integers $n$, the equality
\begin{equation} \label{recursive}
s_{n+1}=\frac{2n-1}{2n+1}\cdot s_n + \frac{4n+4}{3}.
\end{equation}
holds.
\end{lemma}

\begin{proof}
In order to prove (\ref{recursive}), we need the following simple
facts. 

\begin{proposition} \label{facts}
\begin{enumerate}
\item For all positive integers $n$, and all indices $i\neq j$ that
satisfy
$1\leq i,j\leq n$,  
the equality \[E(Y_{n+1,i}Y_{n+1,j})=\frac{2n-1}{2n+1}E(Y_{n,i}Y_{n,j})\]
holds.
\item For all positive integers $n$ and all indices $1\leq i\leq n$, 
the equality \[E(Y_{n+1,i})=\frac{2n}{2n+1}E(Y_{n,i})\] holds.
\item For all indices $i\leq n+1$, the equality
\[E(Y_{n+1,i}Y_{n+1,n+1})=E(Y_{n+1,i})\] holds. In particular, 
$E(Y_{n+1,n+1})=1$. 
\end{enumerate}
\end{proposition}

\begin{proof}
\begin{enumerate}
\item In order to get an element of $Q_{n+1}$ in which $i$ and $j$ are both
plateaux, take an element of $Q_n$ in which  $i$ and $j$ are both plateaux,
and insert two consecutive copies of $n+1$ into any of the $2n-1$ available
places, that is, anywhere but between the two copies of $i$ or the two 
copies of $j$.
\item In order to get an element of $Q_{n+1}$ in which $i$ is a plateau,
insert  two consecutive copies of $n+1$ into any of the $2n$ available
slots, that is, anywhere  but between the two copies of $i$.
\item Obvious since $n+1$ is always a plateau in elements of $Q_{n+1}$.
\end{enumerate} 
\end{proof}

We return to proving Lemma \ref{recu}.  

Note that 
$ s_{n+1}= \sum_{1\leq i,j\leq n+1} 
E(Y_{n+1,i}Y_{n+1,j})$. The latter can be split into partial sums 
based on whether $i$ or $j$ are equal to $n+1$ as follows. 
 \[ s_{n+1}= \sum_{1\leq j\leq n+1} E(Y_{n+1,n+1}Y_{n+1,j})+
 \sum_{1\leq i \leq n} E(Y_{n+1,i}Y_{n+1,n+1})\]
\[+\sum_{1\leq i,j\leq n} 
E(Y_{n+1,i}Y_{n+1,j}) .\]
Based on part 3 of Proposition \ref{facts}, this simplifies to
 \[s_{n+1}=  \sum_{1\leq j\leq n+1} E(Y_{n+1,j}) +  
\sum_{1\leq i \leq n} E(Y_{n+1,i})+\sum_{1\leq i,j\leq n \atop i\neq j} 
E(Y_{n+1,i}Y_{n+1,j}) \] 
\[+\sum_{1\leq i\leq n}E(Y_{n+1,i}).
\]
Now note that the first sum on the right-hand side is just $E(Y_{n+1})$,
the second sum is $E(Y_{n+1}-Y_{n+1,n+1})=E(Y_{n+1})-1$, use part 
1 of Proposition \ref{facts} on the third sum, and part 2 of 
 Proposition \ref{facts} on the fourth sum to get
\[s_{n+1}=2E(Y_{n})-1+ \frac{2n-1}{2n+1}\left(s_n-E(Y_n)\right)
+\frac{2n}{2n+1}E(Y_n).\]
Recalling from Corollary \ref{expect} that $E(Y_n)=\frac{2n+1}{3}$,
this reduces to (\ref{recursive}).
\end{proof}

Using the recursive formula proved in  Lemma \ref{recu}, it is routine
to prove that
\begin{equation}\label{square}
s_n=E(Y_n^2)=\frac{8n^3+6n^2-2n-3}{18n-9}.\end{equation}
Therefore, $\Var(Y_n)=s_n-E(Y_n)^2=\frac{2n^2-2}{18n-9}$ as claimed. 
\end{proof}

\begin{theorem} 
The distribution of the number of plateaux of elements of $Q_n$ converges
to a normal distribution as $n$ goes to infinity. That is, 
$\tilde{Y}_n\rightarrow N(0,1)$. 
\end{theorem}

\begin{proof} Let $X_n=Y_n$, 
and let $g_n(x)=\frac{1}{(2n-1)!!}C_n(x)$. Then
 Theorem \ref{realzeros} and Theorem \ref{varform}
 show that the conditions of Theorem \ref{bender} are 
satisfied, and the claim follows from Theorem \ref{bender}.
\end{proof}

\section{Remarks}

Corollary \ref{expect} shows 
that $E(Y_n)=(2n+1)/3$. It is not difficult to prove that
$E(Y_{n,n-i})=\prod_{j=1}^i \frac{2n-2j}{2n-2j+1}$
By the linearity of expectation this proves
the interesting identity 
\[\sum_{i=0}^{n-1} \prod_{j=1}^i\frac{2n-2j}{2n-2j+1} = \frac{2n+1}{3},\]
where the empty product (indexed by $i=0$) is considered to be 1.

The proof of the equidistribution of the descent and plateau statistics we
gave is very simple, but it is of recursive nature. A direct bijective
proof  has recently been given by  Hyeong-Kwan  Ju \cite{ju}. 
\vskip 1 cm
{\bf \centerline{Acknowledgment}}
I am indebted to Svante Janson, who pointed out an error in an earlier 
version of this paper, which led to an improvement of my results. 
 I am grateful to Ira Gessel for having taken the time to show me some
earlier unpublished work on the subject.

\end{document}